%
%
\BlackBoxes
\documentstyle{amsppt}
\magnification=1200
\topmatter
\title
Upward extension of the Jacobi matrix for orthogonal polynomials
\endtitle
\author
Andr\'e Ronveaux and Walter Van Assche
\endauthor
\affil
Facult\'es Universitaires Notre Dame de la Paix, Namur \\
Katholieke Universiteit Leuven
\endaffil
\address
D\'epartement de Physique, Facult\'es Universitaires N.D.P.,
Rue de Bruxelles 61, B-5000 Namur (BELGIUM)
\endaddress
\email
aronveaux\@cc.fundp.ac.be
\endemail
\address
Department of Mathematics, Katholieke Universiteit Leuven,
Celestijnenlaan 200\,B, B-3001 Heverlee (BELGIUM)
\endaddress
\email
Walter.VanAssche\@wis.kuleuven.ac.be
\endemail
\thanks
The second author is a Senior Research Associate of the Belgian National Fund
for Scientific Research
\endthanks
\date March 3, 1995 \enddate
\keywords
Orthogonal polynomials, Jacobi matrices, perturbations
\endkeywords
\subjclass
42C05
\endsubjclass
\rightheadtext{Upward extension of a Jacobi matrix}

\abstract
Orthogonal polynomials on the real line always satisfy a three-term recurrence
relation. The recurrence coefficients determine a tridiagonal semi-infinite
matrix
(Jacobi matrix) which uniquely characterizes the orthogonal polynomials. We
investigate new orthogonal polynomials by adding to the Jacobi matrix $r$
new rows and columns, so that the original Jacobi matrix is shifted downward.
The $r$ new rows and columns contain $2r$ new parameters  and the
newly obtained orthogonal polynomials thus correspond to an upward extension of
the Jacobi matrix. We give an explicit expression of the new orthogonal
polynomials in terms of the original orthogonal polynomials, their
associated polynomials and the $2r$ new parameters, and we give a fourth
order differential equation for these new polynomials when the original
orthogonal polynomials are classical. Furthermore we show
how the orthogonalizing measure for these new orthogonal polynomials can be
obtained and work out the details for a  one-parameter family
of Jacobi polynomials for which the associated polynomials are again Jacobi
polynomials.
\endabstract

\endtopmatter

\document
\head 1. Introduction \endhead
The construction of families of orthogonal polynomials on the real line from
a given system of orthogonal polynomials (or from a given orthogonalizing weight
$\mu$) has been the subject of various investigations
\cite{6}, \cite{4}, \cite{14}, \cite{11}, \cite{3}, \cite{15}, \cite{7},
\cite{9}, \cite{16}, \cite{25}.
Let $P_n$ $(n=0,1,2,\ldots)$ be a sequence of monic orthogonal polynomials
on the real line, with orthogonality measure $\mu$, then these polynomials
satisfy a three-term recurrence relation
$$  P_{n+1}(x) = (x-b_n) P_n(x) - a_n^2 P_{n-1}(x),
  \qquad n \geq 0,  \tag 1.1 $$
with $b_n \in {\Bbb R}$ and $a_n^2 > 0$ and initial conditions
$P_0=1$, $P_{-1}=0$. The corresponding orthonormal polynomials are
$$   p_n(x) = \frac{1}{a_1 a_2 \ldots a_n} P_n(x), $$
and they satisfy the three-term recurrence relation
$$  x p_{n}(x) = a_{n+1} p_{n+1}(x) +b_n p_n(x) + a_n p_{n-1}(x),
  \qquad n \geq 0. \tag 1.2 $$
Putting the recurrence coefficients in an infinite tridiagonal matrix gives
the Jacobi matrix
$$  J =
   \pmatrix
    b_0   & 1     &        &         & \\
    a_1^2 & b_1   & 1      &         & \\
          & a_2^2 & b_2    & 1       & \\
          &       & \ddots & \ddots  & \ddots
    \endpmatrix  , $$
or the symmetric Jacobi matrix
$$  J_s =
   \pmatrix
    b_0   & a_1   &        &         & \\
    a_1   & b_1   & a_2     &         & \\
          & a_2   & b_2    &  a_3     & \\
          &       & \ddots & \ddots  & \ddots
    \endpmatrix  . $$
Interesting new Jacobi matrices can be obtained by deleting the first $r$
rows and columns, and the corresponding orthogonal polynomials are then
the associated polynomials of order $r$, denoted by $p_n^{(r)}$.
Instead of deleting rows and columns we will add $r$ new rows and columns
at the beginning of the Jacobi matrix, thus introducing $2r$ new parameters.
We call the corresponding new orthogonal polynomials anti-associated of order
$r$ and denote them by $p_n^{(-r)}$. In Section 4 we will give explicit
formulas for these anti-associated polynomials in terms of the original
orthogonal polynomials and the associated polynomials. We also give a
fourth order differential equation in case the original orthogonal polynomials
are classical. The construction of the orthogonality measure for the
anti-associated polynomials is done in Section 5. The analysis requires
the knowledge of the orthogonality measure of the original system and
of the associated polynomials, and also the asymptotic behavior of
both $P_n$ and $P_n^{(1)}$ as $n\to \infty$. These things are known
for a particular one-parameter family of Jacobi polynomials which we
will call Grosjean polynomials (after Grosjean who gave a number of
interesting properties of this family in \cite{9}) and which are
considered in Section 3.

The construction of new orthogonal polynomials by changing and shifting
recurrence coefficients has been under investigation by others.
If the original system of orthogonal polynomials has constant recurrence
coefficients, then one essentially deals with Chebyshev polynomials of the
second kind. Changing a finite number of the recurrence coefficients
leads to Bernstein-Szeg\H{o} weights \cite{22, \S 2.6}, i.e., weights of the
form $\sqrt{1-x^2}/\rho(x)$, where $\rho$ is a positive polynomial on the
orthogonality interval $[-1,1]$. Such finite perturbations have been considered
among others by Geronimus \cite{6}, Grosjean \cite{7},
Sansigre and Valent \cite{21}, and Zabelin \cite{25}.
The limiting case, when
the number of changed recurrence coefficients tends to infinity (with changes
becoming smaller) is treated by Geronimo and Case \cite{4}, Dombrowski and
Nevai \cite{3}, and in \cite{24}.
Instead of starting with constant recurrence coefficients
one can also start with periodic recurrence coefficients. Finite
perturbations of periodic recurrence coefficients are considered by Geronimus
\cite{6}, Grosjean \cite{8} and the limiting case by Geronimo and Van
Assche \cite{5}. Changing a finite  number of
recurrence coefficients of a general system of orthogonal polynomials
has been investigated Marcell\'an et al.\ \cite{11} and the limiting case
by Nevai and Van Assche \cite{15}. Associated polynomials
correspond to deleting rows and columns of the Jacobi matrix, or equivalent
to a positive shift in the recurrence coefficients. This situation is rather
well known since it corresponds to numerator polynomials in Pad\'e
approximation. A general study of such a shift in the recurrence coefficients
can be found in Bel\-meh\-di \cite{1} or Van Assche \cite{23}, and for
Jacobi polynomials we refer to Grosjean \cite{9} and Lewanowicz \cite{10}.
The case where
both a shift in the recurrence coefficients is made together with
changing a finite number of coefficients is treated by Nevai  \cite{14}
and Peherstorfer \cite{16}. In particular Nevai shows that when the original
system is orthogonal on an interval $\Delta$ with weight function $w$, then
the orthogonal polynomials with weight function $w(x)/|f_+(x)|^2$ on
$\Delta$, where $f(z)=a + S(Bw,z)$ and $f_+(x) = \lim_{\epsilon \to 0+}
f(x+i\epsilon)$, with $S(Bw,z)$ the Stieltjes transform of the function
$Bw$ and $B$ a polynomial of degree $\ell$, have recurrence coefficients
which can be obtained by shifting the original recurrence coefficients and
changing a finite number of the initial recurrence coefficients.
Note that this generalizes the Bernstein-Szeg\H{o} polynomials.

\head 2. Background and notation \endhead
Let us consider the family of monic orthogonal polynomials
$P_n$, $n=0,1,2,\ldots$, where $P_n$ has degree $n$, defined by the three-term
recurrence relation
$$  P_{n+1}(x) = (x-b_n) P_n(x) - a_n^2 P_{n-1}(x),
  \qquad n \geq 0, a_n^2 \neq 0, \tag 2.1 $$
with $P_{-1}=0$ and $P_0=1$. The sequences $b_n$ $(n \geq 0)$ and $a_n^2$ $(n
\geq 1)$ also generate the {\it associated} monic polynoimials of order $r$
($r$ a positive integer) $P_n^{(r)}$, $n=0,1,2,\ldots$, by the shifted
recurrence
$$  P_{n+1}^{(r)}(x) = (x-b_{n+r}) P_n^{(r)}(x) - a_{n+r}^2 P_{n-1}^{(r)}(x),
  \qquad n \geq 0, $$
with $P_{-1}^{(r)}=0$ and $P_{0}^{(r)}=1$.
The recurrence relation \thetag{2.1} can also be written in operator form
using the non-symmetric Jacobi matrix $J$ given by
$$  J = \left(
   \matrix
   \matrix
    b_0   & 1     &        &         \\
    a_1^2 & b_1   & 1      &         \\
          & a_2^2 & b_2    & 1       \\
          &       & \ddots & \ddots
    \endmatrix  &
    \matrix
            & & & &\\
            & & & &\\
            & & & &\\
     \ddots & \ \  & \ \ &  &\quad
    \endmatrix   \\
    \matrix
     \ \  &  \ \  & \ \ \  &  a_r^2  \\
          &       &        &         \\
          &       &        &         \\
          &       &        &
    \endmatrix   &  \left|
  \overline{  \matrix
          b_r        & 1       &   & \\
          a_{r+1}^2 & b_{r+1} & 1 & \\
          & \ddots & \ddots & \ddots
     \endmatrix} \right.
        \endmatrix  \right) ,
         $$
as
$$   x \vec{P} = J \vec{P}, $$
where $\vec{P} = (P_0,P_1,P_2,\cdots)^t$.
In the same way the recurrence relation defining the associated polynomials
of order $r$ is represented by
$$   x\vec{P}^{\,(r)} = J^{(r)} \vec{P}^{\,(r)}, $$
where $\vec{P}^{\,(r)}= (P_0^{(r)},P_1^{(r)},P_2^{(r)},\cdots)^t$ and
$J^{(r)}$ is the Jacobi matrix obtained by deleting the first $r$ rows and
columns of $J$.

Observe that $J^{(r)} = J$ for every integer  $r$ if and only if
$b_n=b$ for $n \geq 0$ and $a_n^2 = a^2$ for $n > 0$. The corresponding family
of orthogonal polynomials are, up to standardization of the orthogonality
interval to $[-1,1]$, the monic Chebyshev polynomials of the second kind
$u_n(x)=2^{-n} U_n(x)$, for which $b=0$ and $a^2=1/4$. This fundamental
family is a particular case ($\alpha=\beta=1/2$) of the Jacobi family
$P_n^{(\alpha,\beta)}(x)$, which are orthogonal polynomials with respect to the
weight function $(1-x)^\alpha(1+x)^\beta$ on $[-1,1]$. Closely related to the
family $U_n(x)$ are the Chebyshev polynomials of the first kind $T_n(x)$
$(\alpha=\beta=-1/2)$, of the third kind $V_n(x)$ $(\alpha=1/2, \beta=-1/2)$
and of the fourth kind $W_n(x)$ $(\alpha=-1/2, \beta=1/2)$. These monic
Chebyshev polynomials satisfy the same recurrence relation
$$ P_{n+1}(x) = x P_n(x) - \frac14 P_{n-1}(x), \qquad n  \geq 2, \tag 2.2 $$
but with different initial conditions for $P_1$ and $P_2$. These initial
conditions are hidden in the Jacobi matrices $J_T, J_V, J_W$ for respectively
the polynomials $T_n$, $V_n$ and $W_n$, which are explicitly
$$  J_T = \pmatrix  0 & 1 \\
                    \frac12 & J_U \endpmatrix, \quad
  J_V = \pmatrix  -\frac12 & 1 \\
                    \frac14 & J_U \endpmatrix, \quad
  J_W = \pmatrix  \frac12 & 1 \\
                    \frac14 & J_U \endpmatrix, $$
where $J_U$ is the Jacobi matrix for the Chebyshev polynomials of the second
kind $U_n$. These four families can be embedded in the more general
two-parameter family represented by the Jacobi matrix
$$  J^{(-1)} = \pmatrix b & 1 \\
                        c & J_U \endpmatrix, $$
and by definition the corresponding monic polynomials $P_n$ satisfy the
recurrence relation \thetag{2.2}
with initial conditions $P_1(x)=x-b$, $P_2(x)=xP_1(x)-cP_0(x)$.
All the associated polynomials for $r \geq 1$ are the same:
$$ P_n^{(r)}(x) = u_n(x), \qquad n \geq 0, r \geq 1 . $$

\head 3. Grosjean polynomials \endhead
The families $\{t_n,\ n\geq 0\}$ and $\{u_n,\ n \geq 0\}$ of Chebyshev
polynomials have an interesting extension given by Grosjean. Let us consider
first the monic Jacobi family for which the parameters satisfy
$\alpha+\beta=-1$, and define $G_0^\alpha = 1$,
$$  G_n^\alpha(x) = c_n P_n^{(\alpha,-1-\alpha)}(x), \qquad -1 < \alpha < 0, $$
where $c_n$ is a constant making this a monic polynomial, i.e.,
$c_n = 2^n/{2n-1 \choose n}$. These are monic orthogonal polynomials
with respect to the weight function
$$   w_G(x) = \frac{\sin (-\pi\alpha)}{\pi}
\left( \frac{1-x}{1+x} \right)^\alpha \frac{1}{1+x}, \qquad
    -1 < x < 1. \tag 3.1 $$
A second family of monic Jacobi polynomials with $\alpha+\beta=1$ is denoted by
$$  g_n^\alpha(x) = d_n P_n^{(\alpha,1-\alpha)}(x), \qquad -1 < \alpha < 2, $$
with $d_n = 2^n {2n+1 \choose n}^{-1}$, which are monic orthogonal
polynomials with weight function
$$   w_g(x) = \frac{\sin \pi \alpha}{\alpha(1-\alpha)\pi}
\left( \frac{1-x}{1+x} \right)^\alpha (1+x), \qquad
    -1 < x < 1. \tag 3.2 $$
Obviously we have the special cases $G_n^{-1/2}(x) = t_n(x)$ and
$g_n^{1/2}(x) = u_n(x)$. We will refer to the polynomials $G_n^\alpha$
and $g_n^\alpha$ as the {\it Grosjean polynomials} of the first and second
kind respectively. This terminology is justified since Grosjean \cite{9}
showed, using
direct verification, that
$$   (G_n^\alpha)^{(1)} = g_n^{-\alpha} , \tag 3.3 $$
and Ronveaux \cite{19} gave the generalization of the differential link
$t_n'(x) = nu_{n-1}(x)$ between the Chebyshev polynomials
$$  (G_n^\alpha)'(x) = (-1)^{n-1} n g_{n-1}^{-\alpha}(-x) =
     n g_{n-1}^{1+\alpha}.  \tag 3.4 $$
Recall that for $\alpha=-1/2$ we have $u_n(-x) = (-1)^n u_n(x)$.
Grosjean actually shows \cite{9, p.~275} that the only Jacobi polynomials
for which the associated polynomials are again Jacobi polynomials
are the Chebyshev polynomials of the first, second, third and fourth kind,
and the Grosjean polynomials.
Property \thetag{3.3} can easily be checked from the recurrence coefficients,
which are
$$ b_n = \frac{2\alpha+1}{4n^2-1}, \quad a_n^2 =
\frac{(n+\alpha)(n-1-\alpha)}{(2n-1)^2}, \qquad
\text{for } G_n^\alpha, \tag 3.5a $$
and
$$ b_n = \frac{-2\alpha+1}{4(n+1)^2-1}, \quad a_n^2 =
\frac{(n+\alpha)(n+1-\alpha)}{(2n+1)^2}, \qquad
\text{for } g_n^\alpha, \tag 3.5b $$
(see, e.g., Chihara \cite{2, p.~220})
by changing $n$ to $n+1$ and $\alpha$ to $-\alpha$ in \thetag{3.5a}, which
gives \thetag{3.5b}. Property \thetag{3.4} follows by the differential
property $(\hat{P}_n^{(\alpha,\beta)})'(x) =
n\hat{P}_{n-1}^{(\alpha+1,\beta+1)}(x)$
for the monic Jacobi polynomials $\hat{P}_n^{(\alpha,\beta)}$ and the
symmetry property $p_n^{(\alpha,\beta)}(-x)
=(-1)^n p_n^{(\beta,\alpha)}(x)$ (see, e.g., Szeg\H{o} \cite{22, p.~63}).
It is already interesting to note that for the Grosjean polynomials
the $a_n^2$ are rational function of $n$ consisting of the ratio of two
{\it quadratic} polynomials in $n$ (as in the Legendre case), whereas
in general for the Gegenbauer family one has a ratio of two cubic polynomials
in $n$ and for the Jacobi polynomials one deals with quartic polynomials in the
degree $n$.

The pair of Grosjean polynomials also give an answer to the following
question. Let $D$ denote differentiation $D=\frac{d}{dx}$ and let
$L_2 = \sigma(x) D^2 + \tau(x) D + \lambda_n$ be the second order
(hypergeometric or degenerate hypergeometric) differential operator
for the classical orthogonal polynomials (Jacobi, Laguerre, Hermite,
Bessel), and $\{P_n,\ n \geq 0\}$ be the corresponding family of orthogonal
polynomials, where $P_n$  corresponds with the eigenvalue
$\lambda_n=-n[(n-1)\sigma''+2\tau']/2$, so that $L_2(P_n) = 0$. If
$L_2^*$ is the formal adjoint of $L_2$,
$$  L_2^* = L_2 + 2[\sigma'(x)-\tau(x)] D + \sigma''-\tau', $$
then for which {\it orthogonal} family of polynomials $\{P_n^*,\ n \geq 0\}$
does one have $L_2^*(P_n^*) = 0$? Of course, the Legendre polynomials
($\sigma'=\tau$) solve this problem since then $L_2=L_2^*$.
For Grosjean polynomials the operator $L_2$ is
$$  L_{G,\alpha,n} = (1-x^2) D^2 + (-1-2\alpha-x) D + n^2, \qquad
L_{G,\alpha,n}(G_n^\alpha) = 0, $$
and
$$   L_{g,\alpha,n} = (1-x^2) D^2 + (1-2\alpha-3x) D + n(n+2), \qquad
L_{g,\alpha,n}(g_n^\alpha) = 0, $$
and we have
$$  L_{G,\alpha,n}^* = L_{g,-\alpha,n-1},  $$
Therefore, if $P_n = G_n^\alpha$, then $P_{n}^* =
g_{n-1}^{-\alpha}$. From \thetag{3.4} we also have
$$   D\ L_{G,\alpha,n} =   L_{g,1+\alpha,n-1}\ D   . $$

The relative position of the zeros $x_{1,n}^\alpha < x_{2,n}^\alpha < \cdots <
x_{n,n}^\alpha$ of $G_n^{\alpha}$ compared to the zeros
$x_{j,n} = - \cos \frac{(2j-1)\pi}{2n}$ of $T_n$ is controlled by a classical
comparison theorem due to Markov \cite{22, Theorem 6.12.2}. The ratio
between the Chebyshev weight
$$    w_T(x) = \frac1\pi \frac{1}{\sqrt{1-x^2}}, \qquad -1 < x < 1 , $$
and the weight $w_G$ for the Grosjean polynomials $G_n^\alpha$ given
by \thetag{3.1} is
$$  \frac{w_T(x)}{w_G(x)} = \text{const}\times \left(
     \frac{1+x}{1-x} \right)^{\alpha+1/2} . $$
Now $(1+x)/(1-x)$ is an increasing function of $x$ on the interval $[-1,1]$,
hence for $-1/2 < \alpha < 0$ the ratio $w_T/w_G$ is an
increasing function on $[-1,1]$, and consequently we have the following
inequalities
$$     x_{j,n}^\alpha < - \cos \frac{(2j-1)\pi}{2n}, \qquad j=1,2,\ldots,n,\
   -1/2 < \alpha < 0. $$
For $-1 < \alpha < -1/2$ the ratio $w_T/w_G$ is decreasing, and thus the
inequalities for the zeros are reversed
$$     x_{j,n}^\alpha > - \cos \frac{(2j-1)\pi}{2n}, \qquad j=1,2,\ldots,n,\
   -1 < \alpha < -1/2. $$
Similar conclusions can be made for the zeros $y_{j,n}^\alpha$
$(j=1,2,\ldots,n)$ of Grosjean
polynomials $g_n^\alpha$ of the second kind as compared to the zeros
$-\cos \frac{j\pi}{n+1}$ $(j=1,2,\ldots,n)$
of the Chebyshev polynomials of the second kind.
The ratio of the two weights is
$$  \frac{w_U(x)}{w_g(x)} = \text{const.} \left(
\frac{1+x}{1-x} \right)^{\alpha-1/2}, $$
and thus for $1/2 < \alpha < 2$ this ratio is an increasing function so that
$$     y_{j,n}^\alpha < - \cos \frac{j\pi}{n+1}, \qquad j=1,2,\ldots,n,\
   1/2 < \alpha < 2. $$
For $-1 < \alpha < 1/2$ the inequalities are reversed
$$     y_{j,n}^\alpha > - \cos \frac{j\pi}{n+1}, \qquad j=1,2,\ldots,n,\
   -1 < \alpha < 1/2. $$

Finally we note that when we are dealing with Grosjean polynomials of
the first kind, the product of the weight function of the
orthogonal polynomials and the weight function of the associated polynomials
is constant, as in the case of Chebyshev polynomials of the first kind. This
is not true when dealing with the Chebyshev polynomials of the third and fourth
kind.

\head 4. Anti-associated orthogonal polynomials \endhead
The situation described in \S 2 suggests to construct new families
of orthogonal polynomials, which we will denote by $P_{n+r}^{(-r)}$,
obtained by pushing down a given Jacobi matrix and by introducing in the empty
upper left corner new coefficients $b_{-i}$ $(i=r,r-1,\ldots,1)$ on the
diagonal and new coefficients $a_{-i}^2\neq 0$ $(i=r-1,r-2,\ldots,0)$
on the lower subdiagonal. The new Jacobi matrix is then of the form
$$  J^{(-r)} = \left(
   \matrix
   \matrix
    b_{-r}   & 1     &        &         \\
    a_{-r+1}^2 & b_{-r+1}   & 1      &         \\
          & a_{-r+2}^2 & b_{-r+2}    & 1       \\
          &       & \ddots & \ddots
    \endmatrix  &
    \matrix
            & & & &\\
            & & & &\\
            & & & &\\
     \ddots & \ \  & \ \ &  &\quad
    \endmatrix   \\
    \matrix
     \qquad   &  \qquad  & \qquad  &  \  a_0^2 \\
          &       &        &         \\
          &       &        &         \\
          &       &        &
    \endmatrix   &
    \left|
  \overline{  \matrix
          b_0        & 1       &   & \\
          a_{1}^2 & b_{1} & 1 & \\
          & \ddots & \ddots & \ddots
     \endmatrix} \right.
        \endmatrix  \right) .
         $$
We will call the orthogonal polynomials $P_{n}^{(-r)}$
for this Jacobi matrix $J^{(-r)}$
an\-ti-as\-so\-cia\-ted polynomials for the family $P_n$. They contain $2r$ new
parameters and satisfy
$$   [P_{n+r}^{(-r)}(x)]^{(k)} = P_{n+r}^{(k-r)}(x) . $$
For $r=1$ we have
$$  J^{(-1)} = \pmatrix b_{-1} & 1 \\
                        a_0^{2} & J \endpmatrix, $$
and for $r=2$ we have
$$  J^{(-2)} = \pmatrix b_{-2} & 1 & 0 \\
                        a_{-1}^2 & b_{-1} & 1 \\
                         0   & a_{0}^2 & J  \endpmatrix . $$

This new family of anti-associated polynomials $P_{n}^{(-r)}$
can easily be represented
as a combination of the original family $P_n$ and the associated polynomials
$P_{n-1}^{(1)}$. First, denote by $Q_n$ the orthogonal polynomials
for the finite Jacobi matrix
$$  \pmatrix  b_{-r}     &     1    &   &       &     \\
              a_{-r+1}^2 & b_{-r+1} & 1 &       &     \\
                         & a_{-r+2} &  b_{-r+2} & 1   \\
                         &          & \ddots    & \ddots & \ddots \\
                         &          &           & a_{-1}^2 & b_{-1}
    \endpmatrix , $$
so that they satisfy  $Q_0=1$, $Q_{-1}=0$ and
$$   Q_{n+1}(x) = (x-b_{-r+n}) Q_n(x) - a_{-r+n}^2 Q_{n-1}(x), \qquad
     n \leq r-1. $$
Then, clearly
$$    P_{n}^{(-r)}(x) = Q_n(x), \qquad 0 \leq n \leq r. $$
For $n > r$ the anti-associated polynomials satisfy the three-term
recurrence relation
$$    P_{n+r+1}^{(-r)}(x) = (x-b_{n}) P_{n+r}^{(-r)}(x) - a_{n}^2
P_{n+r-1}^{(-r)}(x), \qquad n \geq 0, $$
so that $P_{n+r}^{(-r)}(x)$ is a solution of the three-term recurrence
relation of the original family $P_n(x)$. The initial conditions
however are $P_r^{(-r)}(x) = Q_r(x)$ and $P_{r-1}^{(-r)}(x) = Q_{r-1}(x)$,
and since every solution of the three-term recurrence relation \thetag{2.1}
is a linear combination of $P_n(x)$ and $P_{n-1}^{(1)}(x)$, we have
$$   P_{n+r}^{(-r)}(x) = A P_n(x) + B P_{n-1}^{(1)}(x), \qquad n \geq 0. $$
Using the initial conditions for $n=0$ and $n=1$ gives
$A=Q_r(x)$ and $B=-a_0^2 Q_{r-1}(x)$, and thus we have
$$   P_{n+r}^{(-r)}(x) = Q_r(x) P_n(x) - a_0^2 Q_{r-1}(x) P_{n-1}^{(1)}(x),
\qquad n \geq 0. \tag 4.1 $$

From this representation it is easy to construct a differential equation
satisfied by the family $P_{n+r}^{(-r)}$ if the original family $P_n$
is itself solution of a linear differential equation of second order,
for instance when $P_n$ are the classical polynomials (Jacobi, Laguerre,
Hermite, Bessel), then they are a solution of the hypergeometric
differential equation
$$  L_2 y \equiv \sigma(x) y'' + \tau(x) y' + \left(
    {-n \over 2} (n-1) \sigma'' - \tau'n \right) y = 0. \tag 4.2 $$
The techniques used in \cite{17}, \cite{18}, \cite{19}, \cite{20} and the
fact that
$$  L_2^* P_{n-1}^{(1)} = (\sigma'' - 2\tau') P_n'  \tag 4.3 $$
easily give a fourth order differential equation satisfied by the
anti-associated polynomials $P_{n+r}^{(-r)}(x)$.
The general technique is the following. Let $B(x) = -a_0^2 Q_{r-1}(x)$ and
put $J(x) = B(x) P_{n-1}^{(1)}(x)$. Then we can  transform the equation
\thetag{4.3} to
$$  L_2^*  \frac{J}{B} = (\sigma'' - 2\tau') P_n',  $$
and introducing the differential operator
$$  \multline
R_2 = \sigma B^2 D^2 + [(2\sigma'-\tau) B^2 - 2\sigma
B B'] D  \\
+\ 2\sigma (B')^2 -\ \sigma B B'' - B
B'(2\sigma'-\tau)-[\frac{\sigma''}{2}(n^2-n-2)+\tau'(1+n)]B^2
  \endmultline $$
this is equivalent to
$$  R_2 J = (\sigma'' - 2\tau') B^3 P_n'. \tag 4.4 $$
From \thetag{4.1} we then have
$$  R_2 P_{n+r}^{(-r)}
   = R_2 [ Q_r P_n ] - (\sigma'' - 2\tau') a_0^6 Q_{r-1}^3 P_n', $$
and eliminating the second derivative using \thetag{4.2} then leads to
$$ R_2 P_{n+r}^{(-r)} = M_0 P_n + N_0 P_n', \tag 4.5 $$
where $M_0$ and $N_0$ are polynomials. Taking the derivative in \thetag{4.5}
and using \thetag{4.2} to eliminate $P_n''$ also gives
$$   \sigma [R_2 P_{n+r}^{(-r)}]' = M_1 P_n + N_1 P_n', $$
where $M_1$ and $N_1$ are polynomials, and repeating this also gives
$$   \sigma [\sigma[ R_2 P_{n+r}^{(-r)}]']' = M_2 P_n + N_2 P_n'. $$
This shows that
$$  \det \pmatrix
  R_2 P_{n+r}^{(-r)} & M_0 & N_0 \\
 \sigma [R_2 P_{n+r}^{(-r)}]' & M_1 & N_1 \\
 \sigma[\sigma [R_2 P_{n+r}^{(-r)}]']' & M_2 & N_2
 \endpmatrix = 0, \tag 4.6 $$
which is the desired fourth order differential equation. When $P_n(x) =
G_n^\alpha(x)$, then $\sigma''-2\tau'=0$ so that \thetag{4.4} simplifies
and becomes homogeneous. The differential equation \thetag{4.6}
however remains one of the fourth order, except when $\alpha=-1/2$, because
then
$$  P_{n+r}^{(-r)}(x) = A_0(x) T_n(x) + B_0(x) T_n'(x), $$
where $A_0$ and $B_0$ are polynomials. Similar as in the above reasoning
we then get a second order differential equation
$$  \det \pmatrix
   P_{n+r}^{(-r)} & A_0 & B_0 \\
 \sigma [P_{n+r}^{(-r)}]' & A_1 & B_1 \\
 \sigma[\sigma [P_{n+r}^{(-r)}]']' & A_2 & B_2
 \endpmatrix = 0,  $$
which is equivalent to the equation given in \cite{21}.

\head 5. Construction of the orthogonality measure \endhead

If we use probability measures throughout the analysis, and if
we use lower case $p$ and $q$
for the orthonormal polynomial, then we
have $p_n = \gamma_n P_n$, where
$$  \gamma_n = (a_1a_2 \cdots a_n)^{-1}. $$
Similarly
$$  q_n(x) = \frac{Q_n(x)}{a_{-r+1}a_{-r+2}\cdots a_{-r+n}}, \qquad n \leq r, $$
and $p_{n-1}^{(1)} = \gamma^{(1)}_n P_{n-1}^{(1)}$ where
$\gamma_n^{(1)} = a_1 \gamma_n$. Thus the orthonormal anti-associated
polynomials are
$$  p_{n+r}^{(-r)}(x) = (a_{-r+1} \cdots a_{0})^{-1} \gamma_n P_{n+r}^{(-r)}(x),
$$
and using this in \thetag{4.1} gives
$$  p_{n+r}^{(-r)}(x) = q_r(x) p_n(x) - \frac{a_0}{a_1}
q_{r-1}(x) p_{n-1}^{(1)}(x). \tag 5.1 $$
The orthonormal polynomials are useful in obtaining the weight function
for the anti-associated polynomials. Indeed, we can compute the weight function
using Christoffel functions
$$  \lambda_n(x) = \left( \sum_{j=0}^{n} p_j^2(x) \right) ^{-1}     $$
by means of the following result of M\'at\'e, Nevai and Totik
\cite{12}

\proclaim{Theorem MNT}
Suppose $p_n(x)$ are orthonormal polynomials for a measure $\mu$ on
$[-1,1]$ and let $\mu$ belong to Szeg\H{o}'s class, i.e.,
$$  \int_{-1}^1 \frac{\log \mu'(x)}{\sqrt{1-x^2}} \, dx > -\infty. $$
Then
$$  \lim_{n\to \infty}  n \lambda_n(x) = \pi \mu'(x) \sqrt{1-x^2} \tag 5.2 $$
holds almost everywhere on $[-1,1]$.
\endproclaim

The Szeg\H{o} condition can be relaxed and in fact it suffices to assume
that $\mu$ is a regular measure on $[-1,1]$ (i.e., $\text{supp}(\mu)
= [-1,1]$ and $\lim_{n \to \infty} \gamma_n^{1/n} = 2$) and
$$  \int_{a}^b \log \mu'(x) \, dx > -\infty $$
in order that \thetag{5.2} holds almost everywhere on $[a,b]
\subset (-1,1)$ \cite{12, Thm.~8}. We will assume these conditions
for the measure $\mu$ and moreover we allow the addition of a finite
number of mass points to $\mu$.  Then \thetag{5.2} will still hold
almost everywhere on $[a,b]$. Indeed, if we add a mass point $c$ to $\mu$
then by the extremum property
$$  \lambda_n(x;\mu) = \min_{q_n(x)=1} \int q_n^2(t) \, d\mu(t), $$
where the minimum is taken over all polynomials $q_n$ of degree at most
$n$ which take the value $1$ at the point $x$, we see that for
the measure $\mu_c=\mu + \epsilon \delta_c$ (i.e., the measure $\mu$
to which we add a mass point at $c$ with mass $\epsilon$)
$$ \align
   \lambda_n(x;\mu_c) &= \min_{q_n(x)=1} \left( \int q_n^2(t) \, d\mu(t)
   + \epsilon q_n^2(c) \right)  \\
   &\geq \min_{q_n(x)=1} \int q_n^2(t) \, d\mu(t) \\
   &= \lambda_n(x;\mu),
  \endalign   $$
so that $\liminf_{n \to \infty} n \lambda_n(x;\mu) \geq \pi \mu'(x)
\sqrt{1-x^2}$, almost everywhere on $[a,b]$. On the other hand, consider
the polynomial $q_n(t) = (t-c) r_{n-1}(t) / (x-c)$, where $r_{n-1}$
is the minimizing polynomial for the measure $d\nu(t) = (t-c)^2\, d\mu(t)$,
then for $x \neq c$
$$  \align
\lambda_n(x;\mu_c) &\leq  \frac{1}{(x-c)^2} \int (t-c)^2 r_{n-1}^2(t) \,
    d\mu(t) \\
      &= (x-c)^{-2} \int  r_{n-1}^2(t) \, d\nu(t) \\
      &= (x-c)^{-2} \lambda_{n-1}(x;\nu)
   \endalign $$
so that
$\limsup_{n \to \infty} n\lambda_n(x;\mu_c) \leq \lim_{n \to \infty} (x-c)^{-2}
n \lambda_n(x;\nu) = \pi \mu'(x) \sqrt{1-x^2}$ almost everywhere on $[a,b]$,
since $\nu$ is regular and satisfies the
Szeg\H{o} condition on $[a,b]$. Combined with the previous inequality
this gives
$$   \lim_{n \to \infty} n \lambda_n(x;\mu_c) = \pi \mu'(x) \sqrt{1-x^2} $$
almost everywhere on $[a,b]$, independent of the mass point $c$. This
procedure can be repeated, so that adding a finite number of mass points
does not change the behavior in \thetag{5.2}.

In the previous discussion we wanted to allow the measure $\mu$ to have
a larger support than $[-1,1]$ (by allowing mass points outside $[-1,1]$)
but in such a way that \thetag{5.2} still holds. Alternative we can use
a weaker result by Nevai \cite{13, Thm.~54, p.~104} in which the support
of $\mu$ is allowed to be $[-1,1] \cup E$, where $E$ contains at most
a denumerable number of points which can only accumulate at $\pm 1$.

\proclaim{Theorem N}
Suppose $\mu \in \text{M}(0,1)$, i.e., the recurrence coefficients $a_n$
and $b_n$ have asymptotic behavior given by
$$   \lim_{n \to \infty} a_n = 1/2, \quad \lim_{n \to \infty} b_n = 0. $$
Then
$$ \limsup_{n \to \infty} n \lambda_n(x;\mu) = \pi \mu'(x) \sqrt{1-x^2} $$
holds for almost every $x \in \text{supp}(\mu)$.
\endproclaim

With this result, if we can compute the limit of $n \lambda_n(x;\mu)$
for $x \in (-1,1)$, then this limit is almost everywhere equal to
$\pi \mu'(x) \sqrt{1-x^2}$, which allows us to compute the weight
$\mu'$ on $(-1,1)$ without having to worry about the fact that
$\text{supp}(\mu)$ may be larger than $[-1,1]$. The existence
of such a limit is however not guaranteed
by Nevai's theorem, contrary to the theorem of M\'at\'e, Nevai and Totik
where the existence of the limit is in the conclusion of the theorem
but where $\text{supp}(\mu) = [-1,1]$ is required.

If we assume that the original family $p_n$ belongs to
a measure $\mu$ in the class M$(0,1)$,
then the anti-associated
polynomials $p_n^{(-r)}$ will also have a measure $\mu^{(-r)}$
which belongs to the class M$(0,1)$, and thus we know that $\mu^{(-r)}$
has support $[-1,1] \cup E$, where
$E$ is at most denumerable with the only accumulation points at $\pm 1$
\cite{24, Thm.~1,  p.~437}.
In fact, there can be at most $2r$ mass points, since the associated
polynomials of order $r$ of the anti-associated polynomials $p_n^{(-r)}$
are again $p_n$ and these have all their zeros inside $[-1,1]$, and the
interlacing property of orthogonal polynomials and associated orthogonal
polynomials shows that adding one row and one column in the Jacobi matrix to
form the anti-associated polynomials $p_n^{(-1)}$ can at most add a mass point
to the left of $-1$ and to the right of $1$. Adding $r$ rows and columns thus
can add at most $r$ mass points to the left of $-1$ and $r$ mass points to the
right of $1$. If the original measure $\mu$ belongs to the Szeg\H{o} class
(which is a subclass of M$(0,1)$), then
the measure $\mu^{(-r)}$ restricted to $[-1,1]$ also belongs to the
Szeg\H{o} class. This is so because if we denote by $\mu^{(-r)}_{[-1,1]}$ the restriction of
$\mu^{(-r)}$ to the interval $[-1,1]$, then
$\lim_{n \to \infty}
\gamma_n(\mu^{(-r)})/\gamma_n(\mu^{(-r)}_{[-1,1]})$ exists and is strictly positive \cite{13,
Thm.~25 on p.~136}
and since $\gamma_n(\mu^{(-r)}) = \gamma_n(\mu) (a_{-r+1} \cdots a_0)^{-1}$
and $\mu$ belongs to the Szeg\H{o} class, which is equivalent with
the statement that
$\lim_{n \to \infty} \gamma_n/2^n$ exists and is strictly positive,
it follows that $\lim_{n \to \infty} \gamma_n(\mu^{(-r)}_{[-1,1]})/2^n$ exists
and is strictly positive, so that $\mu^{(-r)}_{[-1,1]}$ belongs to the Szeg\H{o} class.

In order to use the result in \thetag{5.2} we observe that
for $n \geq r$
$$ \sum_{j=0}^{n} [ p_j^{(-r)}(x)]^2
  = \sum_{j=0}^{r-1} q_j^2(x) + \sum_{j=0}^{n-r} [p_{j+r}^{(-r)}(x)]^2 , $$
and using \thetag{5.1} this gives
$$ \align
\sum_{j=0}^{n-r} [p_{j+r}^{(-r)}(x)]^2
   &= q_r^2(x) \sum_{k=0}^{n-r} p_k^2(x)
     + \left( \frac{a_0}{a_1} \right)^2 q_{r-1}^2(x) \sum_{k=0}^{n-r-1}
     [p_{k}^{(1)}(x)]^2  \\
   &\quad  - \frac{2a_0}{a_1} q_r(x)q_{r-1}(x) \sum_{k=1}^{n-r}
     p_k(x)p_{k-1}^{(1)}(x).
\endalign $$
In order to be able to compute these sums, we will need to be able
to compute Christoffel functions of associated polynomials and
sums of mixed form containing the product $p_n(x) p_{n-1}^{(1)}(x)$.
This is in general not so easy. However, when we use Grosjean polynomials
of the first kind, then the associated polynomials are Grosjean polynomials
of the second kind, so that we are always dealing with Jacobi
polynomials, for which we can compute these sums, at least as
$n$ becomes large. For the anti-associated polynomials corresponding
to Grosjean polynomials we thus can prove the following result.

\proclaim{Theorem 1}
If the original system of orthogonal polynomials consists of Grosjean
polynomials of the first kind, then the anti-associated polynomials
$p_n^{(-r)}$ are orthogonal with respect to a measure $\mu^{(-r)}$ which
is absolutely continuous on $[-1,1]$ with density
$$   w_r(x) = \frac{\sin(-\pi\alpha)}{\pi}
   \frac{(1-x)^{\alpha}}{(1+x)^{\alpha+1}}
    \left| q_r(x) - a_0 e^{i\alpha\pi} q_{r-1}(x)
\frac{(1-x)^\alpha}{(1+x)^{\alpha+1}} \right|^{-2}. $$
In addition, there may be at most $2r$ mass point outside $[-1,1]$
($r$ to the left of $-1$ and $r$ to the right of\/ $1$). These mass points
are the roots of the equation
$$   q_r(x) - a_0 q_{r-1}(x) \text{ sign } x \
\frac{|x-1|^\alpha}{|x+1|^{\alpha+1}} = 0. $$
\endproclaim

\demo{Proof}
If the original family consists of Grosjean
polynomials, then by \thetag{4.3} and \thetag{3.1} we have for $-1 < x < 1$
$$  \lim_{n \to \infty} \frac1n \sum_{k=0}^{n-r} p_k^2(x)
    =  \frac{1}{\sqrt{1-x^2}}  \frac{1}{\sin (-\pi \alpha)}
    \frac{(1+x)^{\alpha+1}}{(1-x)^\alpha},  \tag 5.3 $$
and this even holds uniformly on closed  subintervals of $(-1,1)$, and
similarly  by using \thetag{3.2} (with $\alpha$ replaced by $-\alpha$)
$$ \lim_{n \to \infty} \frac1n \sum_{k=0}^{n-r-1}
     [p_{k}^{(1)}(x)]^2 =  \frac{1}{\sqrt{1-x^2}}
     \frac{-2\alpha(1+\alpha)}{\sin (-\pi\alpha)}
     \frac{(1-x)^\alpha}{(1+x)^{\alpha+1}} .  \tag 5.4  $$
For the remaining sum, which consists of a mixture of the original polynomials
and the associated polynomials, we can use Darboux's extension of the
Laplace-Heine formula for Legendre polynomials \cite{22, Thm.\ 8.21.8 on
p.\ 196}
$$ \align
\sqrt{n\pi} P_{n}^{(\alpha,\beta)}(\cos \theta) &=
    \left( \sin \frac{\theta}2 \right)^{-\alpha-1/2}
    \left( \cos \frac{\theta}2 \right)^{-\beta-1/2}   \\
  &\quad  \cos\left( [n+(\alpha+\beta+1)/2]\theta -
    \frac{\alpha+1/2}{2} \pi \right) + O(1/n)  ,  \tag 5.5
 \endalign   $$
which holds uniformly for $x=\cos \theta$ on closed intervals of $(-1,1)$.
For Grosjean polynomials $P_n(x) = G_n^\alpha(x)$ we have
$$    p_n(x) =[1+ O(1/n)]\sqrt{\frac{2n\pi}{\sin(-\pi\alpha)}}
  P_{n}^{(\alpha,-1-\alpha)}(x), $$
$$ p_{n-1}^{(1)}(x) = [1+O(1/n)]\sqrt{\frac{-n\pi
\alpha(1+\alpha)}{\sin(-\pi\alpha)}} P_{n-1}^{(-\alpha,1+\alpha)}(x), $$
and hence for $x=\cos \theta$ with $\epsilon \leq \theta \leq \pi-\epsilon$
$$ \multline
p_n(x) = \sqrt{\frac{2}{\sin(-\pi\alpha)}}
    \left( \sin \frac{\theta}2 \right)^{-\alpha-1/2}
    \left( \cos \frac{\theta}2 \right)^{\alpha+1/2}
    \cos \left( n\theta - \frac{\alpha+1/2}{2} \pi \right) \\
    + O(1/n)
\endmultline    $$
and
$$ \multline
p_{n-1}^{(1)}(x) = \sqrt{\frac{-\alpha(1+\alpha)}{\sin(-\pi\alpha)}}
    \left( \sin \frac{\theta}2 \right)^{\alpha-1/2}
    \left( \cos \frac{\theta}2 \right)^{-\alpha-3/2}
    \cos \left( n\theta + \frac{\alpha-1/2}{2} \pi \right) \\
    + O(1/n).
\endmultline    $$
Using this gives
$$ \multline
\lim_{n \to \infty} \frac1n \sum_{k=1}^{n-r}
     p_k(x)p_{k-1}^{(1)}(x)  =
    \frac{\sqrt{-2\alpha(1+\alpha)}}{\sin(-\pi\alpha)}
    \frac{1}{\sin \theta/2 \cos \theta/2} \\
     \times  \lim_{n \to \infty} \frac 1n \sum_{k=1}^{n-r}
    \cos \left(  k\theta - \frac{\alpha+1/2}2 \pi \right)
    \cos \left(  k\theta - \frac{-\alpha+1/2}2 \pi \right)  .
    \endmultline  $$
Using $2 \cos a \cos b = \cos(a+b) + \cos(a-b)$ gives
$$ \align
\frac 1n \sum_{k=1}^{n-r}
    \cos \left(  k\theta - \frac{\alpha+1/2}2 \pi \right)
  &  \cos \left(  k\theta - \frac{-\alpha+1/2}2 \pi \right) \\
   &= \frac{1}{2n} \sum_{k=1}^{n-r}
      \left[ \cos \left( 2k\theta - \pi/2 \right) + \cos \alpha\pi \right]
\endalign      $$
and thus uniformly for $x$ on closed subsets of $(-1,1)$ we have
$$  \lim_{n \to \infty} \frac1n \sum_{k=1}^{n-r}
     p_k(x)p_{k-1}^{(1)}(x)
   = \frac{\sqrt{-2\alpha(1+\alpha)}\cos\pi\alpha}{\sin \theta\sin(-\pi\alpha)}.
   \tag 5.6 $$
Alternatively, \thetag{5.5} can also be obtained by using a Christoffel-Darboux
type formula \cite{1, corollary 2.12} which for orthonormal polynomials
is
$$   \sum_{k=1}^n p_k(x) p_{k-1}^{(1)}(x)
= a_{n+1} \left[
  p_{n+1}'(x)p_{n-1}^{(1)}(x) - p_n'(x)p_n^{(1)}(x) \right] . $$
When the original system consists of Grosjean polynomials
of the first kind $G_n^\alpha$ one knows by \thetag{3.3}
that the associated
polynomials are Grosjean polynomials of the second kind $g_n^{-\alpha}$
and by \thetag{3.4} we know that the derivative of $G_n^\alpha$ is
$n g_{n-1}^{1+\alpha}$. For the orthonormal polynomials this gives
$$    p_n(x) = [1+O(1/n)] \sqrt{\frac{2n\pi}{\sin(-\pi\alpha)}}
  P_{n}^{(\alpha,-1-\alpha)}(x), $$
$$ p_{n}^{(1)}(x) = [1+O(1/n)] \sqrt{\frac{-n\pi
\alpha(1+\alpha)}{\sin(-\pi\alpha)}} P_{n}^{(-\alpha,1+\alpha)}(x), $$
$$  p_n'(x) = [1+O(1/n)]  \sqrt{\frac{2\pi n}{\sin(-\pi\alpha)}}
    \frac{n}{2}  P_{n-1}^{(1+\alpha,-\alpha)}(x), $$
and hence using \thetag{5.5} gives
$$  \align
\lim_{n \to \infty} \frac1n  & \sum_{k=1}^n  p_k(x)  p_{k-1}^{(1)}(x)    \\
= &
 \frac{\sqrt{-2\alpha(1+\alpha)}}{\sin(-\pi\alpha) \sin^2\theta}
 \lim_{n \to \infty}  \left[
  \cos \left( (n+1)\theta - \frac{\alpha +3/2}2 \pi \right)
  \cos \left( n\theta - \frac{-\alpha +3/2}2 \pi \right)  \right. \\
&\quad  - \left. \cos \left( n\theta - \frac{\alpha +3/2}2 \pi \right)
  \cos \left( (n+1)\theta - \frac{-\alpha +1/2}2 \pi \right)  \right]
  \endalign
$$
Simple trigonometry then gives \thetag{5.6}.
Combining the limiting relations \thetag{5.3}, \thetag{5.4}, and \thetag{5.6}
gives uniformly for $x$ on closed intervals of $(-1,1)$
$$ \align
\lim_{n \to \infty}
     \frac{1}{n} \sum_{j=0}^n  [ p_j^{(-r)}(x)]^2 & =
      \frac{1}{\sqrt{1-x^2}} \frac{1}{\sin(-\pi\alpha)}
      \left[q_r^2(x)
     \frac{(1+x)^{\alpha+1}}{(1-x)^\alpha} \right. \\
  &\qquad    -\ 2 \frac{a_0}{a_1} q_r(x)q_{r-1}(x)
  \cos \alpha \pi \sqrt{-2\alpha(1+\alpha)}  \\
 &\qquad  + \left. \left(\frac{a_0}{a_1} \right)^2 q_{r-1}^2(x)
(-2\alpha)(1+\alpha)
\frac{(1-x)^{\alpha}}{(1+x)^{\alpha+1}} \right]
\endalign $$
Observe now that $a_1 = \sqrt{-2\alpha(1+\alpha)}$, thus we find
$$ \multline \lim_{n \to \infty}
     \frac{1}{n} \sum_{j=0}^n  [ p_j^{(-r)}(x)]^2       \\
  = \frac{1}{\sin(-\pi\alpha)} \frac{1}{\sqrt{1-x^2}}
   \frac{(1+x)^{\alpha+1}}{(1-x)^{\alpha}}
   \left| q_r(x) - a_0 e^{i\alpha\pi} q_{r-1}(x)
\frac{(1-x)^\alpha}{(1+x)^{\alpha+1}} \right|^2 ,
 \endmultline  $$
and this holds uniformly on closed intervals $[a,b] \subset (-1,1)$.
From \thetag{5.2} we can then conclude that the orthogonality
measure $\mu^{(-r)}$ for the anti-associated polynomials has an absolutely
continuous part on $(-1,1)$ with weight function
$$   w_r(x) = \frac{\sin(-\pi\alpha)}{\pi}
   \frac{(1-x)^{\alpha}}{(1+x)^{\alpha+1}}
    \left| q_r(x) - a_0 e^{i\alpha\pi} q_{r-1}(x)
\frac{(1-x)^\alpha}{(1+x)^{\alpha+1}} \right|^{-2}, \tag 5.7 $$
which is the weight function of the original system $p_n(x)$ (Grosjean
polynomials of the first kind) divided by a positive factor containing
the new parameters $b_{-r},\ldots,b_{-1}$ and $a_{-r+1}^2,
\ldots,a_{-1}^2,a_0^2$.
This factor cannot vanish on $(-1,1)$, because then both the real part and the
imaginary part of
$$  q_r(x) - a_0 e^{i\alpha\pi} q_{r-1}(x)
\frac{(1-x)^\alpha}{(1+x)^{\alpha+1}}   $$
need to vanish. The imaginary part can only vanish when $q_{r-1}(x)=0$, and
assuming this, the real part can then only vanish when also $q_r(x) =0$.
This is impossible since two consecutive orthogonal polynomials cannot
have a common zero.

Since the original orthogonal polynomials are Jacobi polynomials, they
will belong to the class M$(0,1)$. Moreover, the original system satisfies
$$  \sum_{k=0}^\infty ( |1-4a_{k+1}^2| + 2 |b_k| ) < \infty, $$
and hence also the new system of anti-associated polynomials satisfies this
{\it trace class} condition. But then we know \cite{24, Theorem 6}
that the orthogonality measure $\mu^{(-r)}$ is absolutely continuous
on $(-1,1)$ and we have obtained the weight function in \thetag{5.7}.
The mass points outside $(-1,1)$ are those points $x \in {\Bbb R} \setminus
(-1,1)$ for which $\sum_{k=0}^\infty [p_k^{(-r)}]^2 < \infty$.
This means that at a mass point $x$ we have $p_n(x) \to 0$, and from
\thetag{5.1} this implies that
$$  \lim_{n \to \infty} p_n(x) \left[ q_r(x) - \frac{a_0}{a_1}
     q_{r-1}(x) \frac{p_{n-1}^{(1)}(x)}{p_n(x)} \right] = 0. $$
For $x \notin [-1,1]$ we know that $|p_n(x)|$ increases exponentially
fast, hence at a mass point we always have
$$   q_r(x) - \frac{a_0}{a_1} q_{r-1}(x) \lim_{n \to \infty}
     \frac{p_{n-1}^{(1)}(x)}{p_n(x)}  = 0 . $$
The limit of the ratio $p_{n-1}^{(1)}(x)/p_n(x)$ can be found
by using Markov's theorem, which states that this limit is the Stieltjes
transform of the measure $\mu$, or it can be obtained from Darboux's
generalization of the Laplace-Heine formula for Legendre polynomials
for $x$ outside $[-1,1]$ \cite{22, Theorem 8.21.7}. Both methods
give
$$   \lim_{n \to \infty} \frac{p_{n-1}^{(1)}(x)}{a_1p_n(x)}
    = \frac{(x-1)^\alpha}{(x+1)^{\alpha+1}}, $$
where the right hand side is to be taken positive if $x > 1$ and negative
if $x < -1$. This means that a mass point satisfies
$$   q_r(x) - a_0 q_{r-1}(x) \frac{(x-1)^\alpha}{(x+1)^{\alpha+1}}
  = 0. $$
Clearly, there can be at most $2r$ mass points, as was indicated earlier.
 \qed
\enddemo

{\bf Remark:}
In order to determine the measure $\mu^{(-r)}$ for the anti-associated polynomials
corresponding to Grosjean polynomials, one can also use the results from
Theorem 3.9 in Peherstorfer \cite{16} using the Cauchy principal value
and the Stieltjes transform. Indeed, the Grosjean polynomials
are one of the seldom cases where a nice explicit expression for the Cauchy
principal value and the Stieltjes transform exists and for this reason the
results of \cite{16} can be applied without problems. Observe that
this technique is basically also the one used by Grosjean in \cite{9}.
Our approach using Christoffel functions has the advantage that it
advoids taking boundary values of a Stieltjes transform or evaluating
a Cauchy principal values and uses only
information of the orthogonal polynomials on the interval $[-1,1]$. In
particular our method also works when appropriate asymptotic information of the
orthogonal polynomials and the associated orthogonal polynomials on
the real line is available.

\head 6. Examples \endhead
The simplest examples occur when we take $\alpha = -1/2$, in which case
the original system consists of Chebyshev polynomials of the first kind.
The weight function for the anti-associated polynomials then becomes
$$  w_r(x) = \frac{1}{\pi} \frac{1}{\sqrt{1-x^2}} \frac{1}{|q_r(x)- ia_0
q_{r-1}/\sqrt{1-x^2}|^2}
=  \frac{1}{\pi} \frac{\sqrt{1-x^2}}{(1-x^2)q_r^2(x) + a_0^2 q_{r-1}^2(x)}  .
$$
Hence this weight function is the weight function of Chebyshev polynomials
of the second kind, divided by a polynomial which is positive on $[-1,1]$.
Such orthogonal polynomials are known as Bernstein-Szeg\H{o} polynomials
\cite{22, \S 2.6}.

In case the original system consists of Chebyshev polynomials of the second
kind, for which all the recurrence coefficients are constant, we have
$p_n = U_n$ and $p_{n-1}^{(1)} = U_{n-1}$. Relation \thetag{5.1}
then becomes
$$  p_{n+r}^{(-r)}(x) = q_r(x) U_n(x) - 2a_0 q_{r-1}(x) U_{n-1}(x). $$
Using $U_n(x) = \sin (n+1)\theta/\sin\theta$ ($x = \cos \theta$),
one easily shows for $-1 < x < 1$
$$ \lim_{n \to \infty} \frac1n \sum_{j=0}^n U_j^2(x) =
    \frac{1}{2 (1-x^2)}  $$
and
$$  \lim_{n \to \infty} \frac1n \sum_{j=1}^n U_j(x) U_{j-1}(x)
   = \frac{x}{2(1-x^2)}, $$
and hence
$$  \frac1n \sum_{j=0}^n [p_{j}^{(r)}(x)]^2
   = \frac1{2(1-x^2)} \left[q_r^2(x) -4a_0 x q_r(x)q_{r-1}(x) + 4a_0^2
     q_{r-1}^2(x) \right], $$
and hence the weight function becomes
$$   w(x) = \frac{2}{\pi} \frac{\sqrt{1-x^2}}{q_r^2(x) -
   4a_0 x q_r(x)q_{r-1}(x) + 4a_0^2 q_{r-1}^2(x)} , $$
i.e., this is again a Bernstein-Szeg\H{o} weight. Observe that it can be written
as
$$   w(x) = \frac{2}{\pi} \frac{\sqrt{1-x^2}}{\left| q_r(x) -
   2a_0 e^{i\theta} q_{r-1}(x)\right|^2}. $$
For the mass points we see that they can occur only for $x \notin [-1,1]$
when
$$ q_r(x) - 2a_0 q_{r-1}(x) \lim_{n\to \infty} \frac{U_{n-1}(x)}{U_n(x)} = 0, $$
and since $U_{n-1}(x)/U_n(x) \to 1/(x+\sqrt{x^2-1})$ for $x \notin [-1,1]$,
where the limit is to be taken positive for $x > 1$ and negative for $x < -1$,
it follows that $x$ is a mass point only when
$$ q_r(x) - 2a_0 q_{r-1}(x) \frac{1}{x+\sqrt{x^2-1}} = 0. \tag 6.1  $$
If we multiply the left hand side by
$q_r(x) - 2a_0 q_{r-1}(x)/(x-\sqrt{x^2-1})$, then this implies
$$   q_r^2(x) - 4a_0 xq_r(x) q_{r-1}(x) + 4a_0^2 q_{r-1}^2(x) = 0. $$
So the mass points are zeros of the polynomial in the denominator of the weight
function, but only those zeros for which \thetag{6.1} holds.

For $r=1$ we can consider the following cases:
\roster
\item $a_0=1/\sqrt2$ and $b_{-1}=0$. In this case the weight function of the
anti-associated polynomials is  $w(x) = 1/(\pi \sqrt{1-x^2})$ and thus we have
the Chebyshev polynomials of the first kind.
\item $a_0=1/2$ and $b_{-1}=-1/2$. Then the weight function is the one for
Jacobi polynomials with $\alpha=1/2$ and $\beta=-1/2$ and thus we have
Chebyshev polynomials of the third kind $V_n(x)$.
\item Finally when $a_0=1/2$ and $b_{-1}=1/2$ the weight function
becomes the one for Jacobi polynomials with $\alpha=-1/2$ and $\beta=1/2$,
so that we have Chebyshev polynomials of the fourth kind $W_n(x)$.
\endroster


\Refs
\widestnumber\no{24}
\ref \no 1
\by S. Belmehdi
\paper On the associated orthogonal polynomials
\jour J. Comput. Appl. Math. \vol 32 \yr 1990 \pages 311--319
\endref

\ref \no 2
\by T. S. Chihara
\book An Introduction to Orthogonal Polynomials
\publ Gordon and Breach \publaddr New York \yr 1978
\endref

\ref \no 3
\by J. Dombrowski and P. Nevai
\paper Orthogonal polynomials, measures and recurrence relations
\jour SIAM J. Math. Anal. \vol 17 \yr 1986 \pages 752--759
\endref

\ref \no 4
\by J. S. Geronimo and K. M. Case
\paper Scattering theory and polynomials orthogonal on the real line
\jour Trans. Amer. Math. Soc. \vol 258 \yr 1980 \pages 467--494
\endref

\ref \no 5
\by J. S. Geronimo and W. Van Assche
\paper Orthogonal polynomials with asymptotically periodic recurrence
coefficients
\jour J. Approx. Theory \vol 46 \yr 1986 \pages 251--283
\endref

\ref \no 6
\by Ya. L. Geronimus
\paper On some finite difference equations and corresponding systems
of orthogonal polynomials
\jour Zap. Mat. Otd. Fiz.-Mat. Fak. i Kharkov Mat. Obsc. \vol 25
{\rm (4)} \year 1957 \pages 87--100 \lang in Russian
\endref

\ref \no 7
\by C. C. Grosjean
\paper The measure induced by orthogonal polynomials satisfying a recursion
formula with either constant or periodic coefficients, part I: constant
coefficients
\jour Acad. Analecta, Koninkl. Acad. Wetensch. Lett. Sch. Kunsten Belgi\"e
\vol 48 {\rm (3)} \yr 1986 \pages 39--60
\endref


\ref \no 8
\by C. C. Grosjean
\paper The measure induced by orthogonal polynomials satisfying a recursion
formula with either constant or periodic coefficients, part II: pure or mixed
periodic coefficients
\jour Acad. Analecta, Koninkl. Acad. Wetensch. Lett. Sch. Kunsten Belgi\"e
\vol 48 {\rm (5)} \yr 1986 \pages 57--94
\endref

\ref \no 9
\by C. C. Grosjean
\paper The weight functions, generating functions and miscellaneous properties
of the sequences of orthogonal polynomials of the second kind
associated with the Jacobi and Gegenbauer polynomials,
\jour J. Comput. Appl. Math. \vol 16  \year 1986 \pages 259--307
\endref

\ref \no 10
\by S. Lewanowicz
\paper Properties of the polynomials associated with the Jacobi polynomials
\jour Math. Comp. \vol 47 \yr 1986 \pages 669--682
\endref

\ref \no 11
\by F. Marcell\'an, J. S. Dehesa and A. Ronveaux
\paper On orthogonal polynomials with perturbed recurrence relations
\jour J. Comput Appl. Math. \vol 30  \yr 1990 \pages 203--212
\endref

\ref \no 12
\by A. M\'at\'e, P. Nevai and V. Totik
\paper Szeg\H{o}'s extremum problem on the unit circle
\jour Ann. of Math. \vol 134 \yr 1991 \pages 433--453
\endref

\ref \no 13
\by P. G. Nevai
\book Orthogonal Polynomials
\bookinfo Memoirs Amer. Math. Soc. {\bf 213}
\publaddr Providence, RI \yr 1979
\endref

\ref \no 14
\by P. Nevai
\paper A new class of orthogonal polynomials
\jour Proc. Amer. Math. Soc. \vol 91  \yr 1984 \pages 409--415
\endref

\ref \no 15
\by P. Nevai and W. Van Assche
\paper Compact perturbations of orthogonal polynomials
\jour Pacific J. Math. \vol 153 \yr 1992 \pages 163--184
\endref

\ref \no 16
\by F. Peherstorfer
\paper Finite perturbations of orthogonal polynomials
\jour J. Comput. Appl. Math. \vol 44 \yr 1992 \pages 275--302
\endref

\ref \no 17
\by A. Ronveaux
\paper Fourth order differential equations for numerator polynomials
\jour J. Phys. A \vol 21 \yr 1988 \pages 749--753
\endref

\ref \no 18
\by A. Ronveaux
\paper 4th order differential equations and orthogonal polynomials of the
Laguerre-Hahn class
\inbook in `Orthogonal Polynomials and their Applications'
\eds C. Brezinski et al.)
\bookinfo IMACS Annals of Computing and Applied Mathematics \vol 9
\publ J. C. Baltzer AG \publaddr Basel \yr 1991 \pages 379--385
\endref

\ref \no 19
\by A. Ronveaux
\paper On a `first integral' relating classical orthogonal polynomials and their
numerator polynomials
\jour Simon Stevin \vol 66 \yr 1992 \pages 159--171
\endref

\ref \no 20
\by A. Ronveaux, A. Zarzo and E. Godoy
\paper Fourth order differential equations satisfied by the generalized
co-recursive of all classical orthogonal polynomials. A study of their
distribution of zeros
\jour J. Comput. Appl. Math. \toappear
\endref

\ref \no 21
\by G. Sansigre and G. Valent
\paper A large family of semi-classical polynomials: the perturbed
Tche\-bi\-chev
\jour J.\ Comput. Appl. Math. \toappear
\endref

\ref \no 22
\by G. Szeg\H{o}
\book Orthogonal Polynomials
\bookinfo Amer. Math. Soc. Colloq. Publ. \vol 23
\publaddr Providence, RI \yr 1975 (4th edition)
\endref

\ref \no 23
\by W. Van Assche
\paper Orthogonal polynomials, associated polynomials and functions of the
second kind
\jour J. Comput. Appl. Math. \vol 37 \yr 1991 \pages 237--249
\endref

\ref \no 24
\by W. Van Assche
\paper Asymptotics for orthogonal polynomials and three-term recurrence
\inbook in `Orthogonal Polynomials: Theory and Practice' \ed P. Nevai
\bookinfo NATO ASI Series C \vol 294
\publ Kluwer \publaddr Dordrecht \yr 1990 \pages 435--462
\endref

\ref \no 25
\by V. V. Zabelin
\paper Polynomial properties
\jour Vestnik Moskov. Univ. Ser. 1, Mat. Meh. nr. 1 \yr 1989  \pages 98--100
\lang in Russian
\transl \jour Mosc. Univ. Math. Bull. \vol 44 \yr 1989 \pages 88--91
\endref





\endRefs
\enddocument